\documentclass[11pt,a4paper,reqno]{amsart}
\usepackage[T1]{fontenc} % delete this line if T1-encoded fonts are missing
\usepackage{a4,amssymb}
\newcommand{\SU}[1][]{\mathrm{SU}_{#1}}
\newcommand{\SUq}{\SU[q]}
\newcommand{\Sq}[1]{\tilde{\mathrm{S}}_q^{#1}}
\let\nob=\nobreakdash
\newcommand{\mbm}[1]{{\let\mathrm=\mathbf\boldmath{#1}}}
\theoremstyle{plain}
\newtheorem{Prop}{Proposition}
\newtheorem{Thm}[Prop]{Theorem}

\newtheorem{Cor}[Prop]{Corollary}
\theoremstyle{definition}

\def\d{\mathrm{d}}
\def\id{\mathrm{id}}
\def\real{\mathrm{I\kern-.18em R}}
\def\nat{\mathrm{I\kern-.18em N}}
\def\DeltaR{\Delta_\mathrm{R}}
\def\PhiR{\Phi_\mathrm{R}}
\newcommand{\cA}{\mathcal{A}}
\newcommand{\Gof}[2][]{{\Gamma^{\otimes#2}_{#1}}}
\newcommand{\Gw}[1][]{\Gamma^\wedge_{\mathrm{#1}}}
\newcommand{\syL}{\sigma_{\mathrm{L}}}
\newcommand{\Gammad}{{(\Gamma,\d)}}
\def\H#1{\Omega_{#1}}
\newcommand{\ESq}[1]{\mathrm{S}_q^{#1}}
\newcommand{\Oq}{\mathrm{O}_q}
\newcommand{\HES}{\H{}}
\newcommand{\RO}[4]{\hat{R}^{#1#2}_{#3#4}}
\newcommand{\ROm}[5][-]{\hat{R}^{#1}{}^{#2#3}_{#4#5}}
\newcommand{\KO}[4]{K^{#1#2}_{#3#4}}
\newcommand{\IO}[4]{I^{#1#2}_{#3#4}}
\newcommand{\CO}[3]{C#1{#2#3}}
\newcommand{\GES}[1]{\Gamma_{#1}}

\def\ylth{0.5pt}\def\yubh{1.5ex}\newdimen\yu\yu1ex\newcommand{\yttv}[9]{%
\raisebox{\yubh}{\smash{\rule[0\yu]{#1\yu}{\ylth}\kern-#1\yu\rule[-1\yu]%
{#1\yu}{\ylth}\kern-#1\yu\rule[-2\yu]{#2\yu}{\ylth}\kern-#2\yu\rule[-3\yu]%
{#3\yu}{\ylth}\kern-#3\yu\rule[-4\yu]{#4\yu}{\ylth}\kern-#4\yu}\rule[-#5\yu]%
{\ylth}{#5\yu}\kern-\ylth\kern1\yu\rule[-#5\yu]{\ylth}{#5\yu}\kern-\ylth%
\kern1\yu\rule[-#6\yu]{\ylth}{#6\yu}\kern-\ylth\kern1\yu\rule[-#7\yu]{\ylth}%
{#7\yu}\kern-\ylth\kern1\yu\rule[-#8\yu]{\ylth}{#8\yu}\kern-\ylth\kern1\yu%
\rule[-#9\yu]{\ylth}{#9\yu}\kern-\ylth\kern-5\yu\kern#1\yu\rule[0\yu]{\ylth}%
{\ylth}\kern-\ylth}}\newcommand{\ytt}[9]{~\yttv#1#2#3#4#5#6#7#8#9~}
\newcommand{\ynull}{\textbf{(0)}}
\allowdisplaybreaks
\begin{document}
\author{Martin Welk}
\address{Martin Welk, Department of Mathematics,
University of Leipzig, Augustusplatz~10, 04109 Leipzig,
Germany --- dr.welk@gmx.net}
\title[Differential Calculus on Quantum Euclidean Spheres]%
{Covariant First Order Differential Calculus on Quantum
Euclidean Spheres}
\begin{abstract}
We study covariant differential calculus on the quantum spheres
$\ESq{N-1}$ which are quantum homogeneous spaces with coactions
of the quantum groups $\Oq(N)$.

The first part of the paper is devoted to first order differential
calculus. A classification result is proved which says that
for $N\ge6$ there exist exactly two
covariant first order differential calculi on $\ESq{N-1}$
which satisfy the classification constraint that the bimodule
of one-forms is generated as a free left module by the differentials
of the generators of $\ESq{N-1}$.
Although the proof given here guarantees completeness of the
classification only for ${N\ge6}$, the calculi themselves exist
for any ${N\ge3}$. Both calculi can also be constructed by a method
introduced by Hermisson. The number of dimensions for
both calculi is by $1$ higher than that of the classical commutative
calculus; in the limit ${q=1}$, one of them can be factorised in such
a way that the commutative calculus is obtained. In the deformed case,
no calculi exist which have the same number of dimensions and
structure of corepresentations as the classical calculus.
In case ${N=3}$, the result is in accordance with the result obtained
by Apel and Schm{\"u}dgen for the Podle{\'s} sphere.

In the second part, higher order differential calculus and symmetry
are treated on the basis of one of the first order calculi.
The relations which hold for two-forms in the universal higher order
calculus extending the underlying first order calculus are given.
A ``braiding'' homomorphism is found that can be used to define
a higher order differential calculus via antisymmetrisation.
The existence of an upper bound for the order of differential forms
is discussed for different choices of higher order calculi.
\end{abstract}
\maketitle
\section{Introduction}
For a number of years, quantum groups as examples of noncommutative
geometric spaces with a rich additional algebraic structure have
received a high interest, and their investigation has made great
progress during the 1990s. As an important prerequisite for an
understanding of their geometrical structure, covariant differential
calculus on quantum groups has been studied to a high level.
Starting from Woronowicz's initiating work~\cite{wor},
construction and classification of covariant first order differential
calculi, and investigation of their properties, were the topic of
many papers, e.g.~\cite{ju},~\cite{ss1},~\cite{ss2}. Higher order
differential calculus was the object of further work, like~\cite{sch},
and basic concept of noncommutative differential geometry were
studied, e.g.~\cite{ac},~\cite{hs}.

In the last few years, investigations concerning covariant differential
calculus on quantum homogeneous spaces---which are the class of
noncommutative spaces next to quantum groups as judged by the mathematical
structure---are increasing; nevertheless, much work has still to be
done until a comparably deep understanding may be reached as for
quantum groups. A general theory even for first order covariant
differential calculus is still lacking. The number of non-trivial
examples with well-understood differential structure is still small.
Methods to construct first order differential calculi on
certain classes of quantum homogeneous spaces are described
in~\cite{schm},~\cite{uh}. A classifications of covariant first order
differential calculi was given for Podle{\'s}' quantum spheres
in~\cite{as}. Earlier work of the author resulted in similar
classification results for the quantum spheres introduced by
Vaksman and Soibelman~\cite{wk1},~\cite{wk2} and for a class of quantum
projective spaces~\cite{wk3},~\cite{wk4}.
The method enabling all of these classifications which will be used
also in this paper is based on screening morphisms of corepresentations
of the underlying quantum group, namely $\SUq(N)$ for all of the
before-mentioned spaces.

The quantum Euclidean spheres $\ESq{N-1}$ which are the subject of
this paper are quantum homogeneous spaces for the quantum groups
$\Oq(N)$. As compared to the $\SUq(N)$\nob-based Vaksman-Soibelman
spheres $\Sq{2N-1}$, this is supposed to be in some sense
the more natural choice in deforming spheres to quantum spheres
but the calculations become much more involved than in the case
of the Vaksman-Soibelman spheres. Nevertheless, it turns out that
in contrast to the great variety of different covariant first order
differential calculi which were found in~\cite{wk1}
for~$\Sq{2N-1}$, the classification theorem~\ref{thfcl}
for the quantum Euclidean spheres
which is the first main result of this paper lists only two
calculi under an appropriate classification constraint.
The outstanding properties of these calculi make the quantum Euclidean
spheres a promising object for further study.
\section{The Quantum Euclidean Spheres \mbm{$\ESq{N-1}$}}
\subsection{Definitions for quantum spaces}
Here and in the following our definitions are in accordance
with~\cite{ks},~\cite{as},~\cite{rtf}.
Let $\cA$ be a coquasitriangular Hopf algebra (This is what
we shall understand by a quantum group throughout the following.)
Let $\cA$ have the comultiplication~$\Delta_{\cA}$ and
the counit~$\varepsilon_{\cA}$.
A pair ${(X,\DeltaR)}$ consisting of an algebra~$X$ and an
algebra homomorphism ${\DeltaR:X\to X\otimes\cA}$
is called a \textsl{(right) quantum space for $\cA$} if
the following equalities hold:
${(\DeltaR\otimes\id_{\cA})\DeltaR=}{(\id_X\otimes\Delta_{\cA})\DeltaR}$;
${(\id_X\otimes\varepsilon_{\cA})\DeltaR=\id_X}$. Then,
$\DeltaR$ is called \textsl{(right) coaction} of~$\cA$ on~$X$.
For brevity, we shall also denote the quantum space itself by~$X$.
A quantum space $(X,\DeltaR)$ (or, briefly, $X$) for $\cA$ is called
\textsl{quantum homogeneous space} if there exists an embedding
${\iota_X:X\to\cA}$ such that ${\DeltaR=\Delta_{\cA}\circ\iota_X}$,
i.e.\ $X$ can be identified with a sub-algebra of $\cA$, and its
coaction is obtained then by restricting the comultiplication.
\subsection{The quantum Euclidean spheres}
Throughout the following, the deformation parameter $q$ is a
positive real number different from $1$ unless stated otherwise.%
\footnote{For most considerations, $q$ could even be arbitrary complex,
except $0$ and roots of unity; the restriction to real numbers
is necessary only where $*$\nob-structures are involved.}
Moreover, the dimension parameter $N$ will always be a natural number,
${N\ge3}$. We shall always sum over pairs of upper and lower indices
from $1$ to~$N$.

Define $X$ to be the algebra with $N$ generators ${x_1, x_2,\dots,x_N}$,
and relations (cf.~\cite{rtf})
\begin{align}
\ROm klij x_kx_l&=q^{-1}x_ix_j +
\frac{q^{N-1}-q^{N-3}}{1+q^{N-2}}\KO klij x_kx_l\text,
\label{eqesr1}\\
\CO^kl x_kx_l&=1\text,
\label{eqesr2}
\end{align}
where
${\ROm{}{}{}{}=\RO{}{}{}{}-(q-q^{-1})\IO{}{}{}{}+(q-q^{-1})\KO{}{}{}{}}$
is the inverse of the R\nob-matrix of the quantum group $\Oq(N)$,
\begin{equation*}
\RO klij=q^{(k=l)-(k=l')}(i=l)(j=k)+(q-q^{-1})(i<l)((i=k)(j=l)-\KO klij)
\text,
\end{equation*}
and
\begin{align*}
\KO klij&=\CO^ij\CO_kl;&\CO^ij=\CO_ij&=q^{-\varrho_i}(i=j')\text.
\end{align*}
Here, ${(i=j)}$ and ${(i<j)}$ denote Kronecker and Heaviside symbols,
resp.; by an apostrophe we denote the mapping ${': i\mapsto i'=N+1-i}$.
The constants $\varrho_i$ are given by ${\varrho_i=\frac N2-i}$ for
${i<\frac{N+1}2}$, $0$ for ${i=\frac{N+1}2}$,
${\frac N2-i+1}$ for ${i>\frac{N+1}2}$.

A $*$\nob-structure on $\ESq{N-1}$ is given by
\begin{equation*}
x_i^*=\CO^ijx_j\text.
\end{equation*}

By ${x_i\mapsto(1+q^{N-2})^{-1/2}\bigl(u^1_i+q^{N/2-1}u^N_i\bigr)}$,
cf.~\cite{ks},~\cite{rtf},
$X$ (as a $*$\nob-algebra) is embedded into $\Oq(N)$; the coaction
\begin{equation*}
\DeltaR:X\to X\otimes\Oq(N),\qquad x_k\mapsto x_i\otimes u^i_k
\end{equation*}
obtained by restricting the comultiplication $\Delta$ of $\Oq(N)$ then
makes $X$ into a quantum homogeneous space for $\Oq(N)$ which we
shall call \textsl{Euclidean quantum sphere} and denote by $\ESq{N-1}$.
\section{First order differential calculus on quantum
Euclidean spheres}
\subsection{Basic definitions}
Let $X$ be an algebra. A \textsl{first order differential calculus}
on $X$ is a pair $\Gammad$ of an $X$\nob-bimodule $\Gamma$ and
a linear mapping ${\d: X\to\Gamma}$ which fulfils the Leibniz rule,
${\d(xy)=\d x\cdot y+x\cdot y}$ for all ${x,y\in X}$, such that
${\d X}$ spans $\Gamma$ as a left module,
i.e.\ ${\mathrm{Lin}\{x\d y~\vert~x,y\in X\}}$.
The elements of $\Gamma$ are called one-forms, while $\d$ is
called differentiation map.

If ${(X,\DeltaR)}$ is a quantum space for a quantum group $\cA$,
a first order differential calculus $\Gammad$ over $X$ is
\textsl{covariant} if
the well-defined linear mapping
${\PhiR:\Gamma\to\Gamma\otimes\cA}$ given by
${\PhiR(\d x)=(\d\otimes\id_{\cA})(\DeltaR(x))}$ and
${\PhiR(x\omega y)=\DeltaR(x)\PhiR(\omega)\DeltaR(y)}$
for all ${x,y\in X}$, ${\omega\in\Gamma}$, and
$\PhiR$ satisfies the identities
${(\PhiR\otimes\id_{\cA})\PhiR=(\id_{\Gamma}\otimes\Delta)\PhiR}$,
${(\id_{\Gamma}\otimes\varepsilon_{\cA})\PhiR=\id_{\Gamma}}$.

An important role in the investigation of covariant differential
calculi over quantum spaces is played by (right-) invariant
one-forms, i.e. those ${\omega\in\Gamma}$ for which
${\PhiR(\omega)=\omega\otimes{\mathbf{1}_{\cA}}}$ is fulfilled.

If for a covariant first order differential calculus $\Gamma$
over $X$ an invariant one-form $\omega_0$ exist such that
${\d x=\omega_0x-x\omega_0}$ holds for all ${x\in X}$, the calculus
$\Gamma$ is called \textsl{inner}.

Finally, if $\Gammad$ is a first order differential calculus
over a $*$\nob-algebra $X$ then $\Gammad$ is called a
\textsl{$*$\nob-calculus} if for any sum ${\sum_k x_k\d y_k}$
where ${x_k,y_k\in X}$ which equals zero in $\Gamma$ even
${\sum_k\d(y_k^*)x_k^*}$ vanishes. In this case it makes sense
to define ${*:\Gamma\to\Gamma}$ via ${(\d x)^*:=\d (x^*)}$.
\subsection{Classification theorem}
We shall consider \textsl{free} covariant first order
differential calculi on $\ESq{N-1}$, i.e.\ calculi $\Gammad$
fulfilling the constraint that the differentials ${\d x_i}$,
${i=1,\dots,N}$, of the generators of $\ESq{N-1}$ generate
the bimodule of one-forms \textsl{as a free left module.}
Our main result is the following classification theorem for
these calculi.
\begin{Thm}\label{thfcl}
For ${N\ge3}$, there exist two free covariant differential calculi
$\GES+$, $\GES-$ on $\ESq{N-1}$. If ${N\ge6}$, any free covariant
differential calculus on $\ESq{N-1}$ is either $\GES+$ or $\GES-$.

The bimodule structure of $\GES\pm$ is given by
\begin{align}
\d x_i\cdot x_j&= \pm\ROm klij x_k\d x_l +(\pm q-1)x_i\d x_j\notag\\
&\quad\quad +\frac{q^{N}-q^{N-2}}{1\mp q^{N-1}}\KO klij x_k\d x_l
+\frac{(1\mp q)(1+q^{N-2})}{1\mp q^{N-1}}\CO^kl x_ix_jx_k\d x_l\text,
\label{eqbimx}
\end{align}
where the upper signs are valid for $\GES+$ and the lower ones
for $\GES-$.
\end{Thm}

The proof of the classification theorem
will be given in section~\ref{sclpr}.

We stress that both differential calculi described in the theorem exist
for any ${N\ge3}$; it's the uniqueness statement that has to
be restricted to ${N\ge6}$ for technical reasons that become
clear in section~\ref{sclpr}.

In the case ${N=3}$, the quantum Euclidean sphere
$\ESq{N-1}$ is isomorphic to one of Podle{\'s}' quantum spheres
$\mathrm{S}_{qc}$, namely that with parameter ${c=0}$.
In~\cite{as} where covariant first order differential calculi
for the Podle{\'s} spheres were classified under a more general
classification constraint, it was shown that
on $\mathrm{S}_{q0}$ there exist exactly two free covariant
first order differential calculi%
\footnote{Note that ${c=0}$ is an exceptional case---for most values
of $c$, exactly one such calculus exists.}%
, which is in accordance with our result.
Besides, it follows that our uniqueness statement is in fact true even
in the case ${N=3}$ which is not covered by our proof.

We continue by discussing some properties of the differential calculi.
First, we note that the bimodule structure can also be given in a form
which allows to transform any given one-form into a \textsl{right}-module
expression.
\begin{Cor}\label{corflr}
The bimodule structure of $\GES+$ and $\GES-$ is described by
\begin{align}
x_i\d x_j&= \pm\ROm klij \d x_k\cdot x_l +
(\pm q-1)\d x_i\cdot x_j\notag\\
&\quad\quad +\frac{q^{N}-q^{N-2}}{1\mp q^{N-1}}\KO klij \d x_k\cdot x_l
+\frac{(1\mp q)(1+q^{N-2})}{1\mp q^{N-1}}\CO^kl \d x_k\cdot x_lx_ix_j\text.
\label{eqbimxr}
\end{align}
\end{Cor}
\textsc{Proof:} The equivalence of
\eqref{eqbimx} and \eqref{eqbimxr} is checked by direct calculation.\qed

An important observation is that in both calculi
there exists a 1\nob-dimensional vector space of invariant one-forms,
namely the multiples of
\begin{equation*}
\HES:=\CO^kl x_k\d x_l\text.
\end{equation*}
One obtains from the bimodule structure given in the theorem that
\begin{equation*}
\HES x_i=x_i\HES \pm\frac{q^{-1}(1\mp q)(1\mp q^{N-1})}{1+q^{N-2}}
\d x_i\text,
\end{equation*}
where again the upper and lower signs refer to $\GES+$ and $\GES-$,
resp. With
\begin{equation*}
\HES':=\pm q\frac{1+q^{N-2}}{(1\mp q)(1\mp q^{N-1})}\HES\text,
\end{equation*}
one therefore has ${\d x=\HES'x-x\HES'}$ for all ${x\in\ESq{N-1}}$,
which implies the following corollary.
\begin{Cor}
$\GES+$ and $\GES-$ are inner calculi.
\end{Cor}
Finally, the question is to be answered whether
our differential calculi are compatible
with the $*$\nob-structure of $\ESq{N-1}$.
\begin{Cor}
$\GES+$ and $\GES-$ are $*$\nob-calculi.
\end{Cor}
\textsc{Proof:}
It is sufficient to prove that
\begin{align}
x_j^* \d x_i^*&= \pm\ROm klij \d x_l^*\cdot x_k^*
+(\pm q-1)\d x_j^*\cdot x_i^*\notag\\
&\quad\quad
+\frac{q^{N}-q^{N-2}}{1\mp q^{N-1}}\KO klij \d x_l^*\cdot x_k^*
+\frac{(1\mp q)(1+q^{N-2})}{1\mp q^{N-1}}\CO^kl
\d x_l^*\cdot x_k^*x_j^*x_i^*
\text.\label{eqbimxst}
\end{align}
Using the defining relations of $\ESq{N-1}$ together with
the definition of $\CO^ij$ and the Leibniz rule,
one checks that
\begin{equation*}
\CO^kl\d x_l^*\cdot x_k^*=\CO^kl\CO^ks\CO^lt\d x_t\cdot x_s
=\CO^ts\d x_t\cdot x_s
\end{equation*}
Moreover, one has ${\ROm klij\CO^ks\CO^lt=\CO^ik\CO^jl\ROm tslk}$
such that \eqref{eqbimxst} becomes
\begin{align*}
\CO^jt\CO^isx_t \d x_s&= \CO^jt\CO^is\Biggl(
\pm\ROm uvts \d x_u\cdot x_v+(\pm q-1)\d x_t\cdot x_s\\*&
+\frac{q^{N}-q^{N-2}}{1\mp q^{N-1}}\KO uvts \d x_u\cdot x_v
+\frac{(1\mp q)(1+q^{N-2})}{1\mp q^{N-1}}\CO^uv
\d x_u\cdot x_vx_tx_s\Biggr)\text,
\end{align*}
which is equivalent to \eqref{eqbimxr}.\qed
\subsection{The classical limit}
Note that the 1\nob-dimensional space of invariant one-forms is still
present in the ``classical limit'', i.e.\ the limit ${q=1}$, of both
calculi while no invariant one-forms exist in the classical commutative
differential calculus. Indeed, the bimodule structure of $\GES+$ in the
limit ${q=1}$ takes the form
\begin{align*}
\d x_i \cdot x_j&=x_j\d x_i-\frac2{N-1}(i=j')\HES+\frac2{N-1}x_ix_j\HES\\
\intertext{while that of $\GES-$ becomes}\d x_i \cdot x_j&=
-x_j\d x_i -2x_i\d x_j+2x_ix_j\HES\text,
\end{align*}
both being noncommutative.

However, there is an important asymmetry between $\GES+$ and $\GES-$
indicating that $\GES+$ stands in a closer relation to the
classical commutative calculus than $\GES-$.
Namely, for $\GES+$ in the classical limit one has
${\HES x=x\HES}$ for all ${x\in\ESq{N-1}}$, such that
$\GES+$ loses in the limit the property of being an inner
calculus; it can be factorised by the additional relation ${\HES=0}$.
The resulting differential calculus is obviously the classical
commutative one, ${\d x_i\cdot x_j=x_j\d x_i}$.
Such a factorisation is not possible for $\GES-$ in the limit ${q=1}$
since we have for $\GES-$ with ${q=1}$ the equality
${\HES x=x\HES-2\d x}$ for all ${x\in\ESq{N-1}}$.%
\footnote{Instead, $\GES-$ admits factorisation over the relation
${\HES=0}$ in the case ${q=-1}$.}

In case ${N\ge6}$,
$\GES+$ is in some sense the best one can get as an approximation
of the classical commutative differential calculus on $\ESq{N-1}$
in the deformed case since there is no covariant first order
differential calculus whose bimodule of one-forms decomposes in
the same way into invariant subspaces for the coaction of $\Oq(N)$
as in the classical case. Namely, this would require (opposite to
the case of a free calculus) that $\HES$ vanishes.

\begin{Cor}\label{crffcl}
Let ${q\in\real\setminus\{0,\pm1\}}$, and ${N\ge6}$. Then there exists no
covariant first order differential calculus $\Gammad$ on $\ESq{N-1}$
with $\Gamma$ generated as a left module by ${\d x_i}$, ${i=1,\dots,N}$,
for which all left-module relations within $\Gamma$ are generated by
the relation ${\CO^ijx_i\d x_j=0}$.
\end{Cor}
The proof of this corollary will be given along with that of
the classification theorem in~\ref{sclpr}.
\subsection{A new basis for one-forms}
From equations \eqref{eqbimx} it is obvious that the
transformation of one-forms into their
left-module expressions is non-linear in $\GES{\pm}$.
For example, in transforming ${\d x_i\cdot x_j}$ to its
left-module expression, besides members of the type ${x_k\d x_l}$
also expressions ${x_kx_lx_m\d x_n}$, with algebra elements
of third degree, are obtained.
Of course, this behaviour of the transformations is a great
obstacle in doing more advanced calculations.
Therefore, the question rises whether the equations describing
the bimodule structure
of the free covariant first order differential calculi on
$\ESq{N-1}$ could be simplified to a purely linear structure
like ${\gamma_i x_j=A^{kl}_{ij}x_k\gamma_l}$
by choosing a different basis $\{\gamma_i\}$
for the bimodule of one-forms.

For reasons of covariance, it is sufficient to coonsider
sets of one-forms of type
${\gamma_i=\beta\d x_i+\alpha x_i\HES}$, with $\alpha$, $\beta$
independent on $i$. Since $\GES{\pm}$
are inner calculi, $\beta$ can't vanish such that we can
restrict ourselves to the ansatz ${\gamma_i=\d x_i+\alpha x_i\HES}$.
\begin{Prop}\label{prga}
Let $\Gammad$ be one of the differential calculi
${(\GES+,\d)}$ and ${(\GES-,\d)}$. Then
there are exactly two real numbers $\alpha$, for which
the set of one-forms
${\{\gamma_i:=\d x_i+\alpha x_i\HES~\vert~i=1,\dots,N\}}$
forms a basis of $\Gamma$ such that
the bimodule structure of $\Gamma$ is described w.r.t.\
the basis $\{\gamma_i\}$ by
\begin{equation*}
\gamma_i x_j=A^{kl}_{ij}x_k\gamma_l
\end{equation*}
with appropriate real coefficients $A^{kl}_{ij}$.

These values are
\begin{align*}
\alpha^+&=-\frac{1+q^{N-2}}{1\mp q^{N-1}}\text,&
\alpha^-&=\pm q\frac{1+q^{N-2}}{1\mp q^{N-1}}
\end{align*}
where the upper signs are valid for ${\Gamma=\GES+}$
while the lower ones refer to ${\Gamma=\GES-}$.
W.r.t.\ the two bases
${\gamma^+_i:=\d x_i+\alpha^+x_i\HES}$,
${\gamma^-_i:=\d x_i+\alpha^-x_i\HES}$
the bimodule structure of $\Gamma$ can be written as
\begin{align}
\gamma^+_ix_j&=\RO  klij x_k\gamma^+_l\text,&
\gamma^-_ix_j&=\ROm klij x_k\gamma^-_l\text.
\label{eqbimg}
\end{align}
\end{Prop}
\textsc{Proof:}
Let ${\Gamma=\GES+}$, and ${\gamma_i:=\d x_i+\alpha x_i\HES}$, with
arbitrary ${\alpha\in\real}$. Because of
\begin{equation*}
\CO^ij x_i\gamma_j=\CO^ij(x_i\d x_j+\alpha\CO^klx_ix_jx_k\d x_l)
=(\alpha+1)\CO^ijx_i\d x_j
\end{equation*}
the inverse substitution is given by
${\d x_i=\gamma_i-(\alpha/(\alpha+1))x_i\CO^klx_k\gamma_l}$.
We calculate
\begin{align}
\gamma_i x_j&=\d x_i\cdot x_j + \alpha x_i\HES x_j\notag\\*
&=\ROm klij x_k\d x_l +
(q-1)\left(1-q^{-1}\alpha\frac{1-q^{N-1}}{1+q^{N-2}}\right)x_i\d x_j
\notag\\*&\quad +
\left(\frac{(1-q)(1+q^{N-2})}{1-q^{N-1}}+\alpha\right)x_ix_j\HES
\frac{q^{N-2}(q^2-1)}{1-q^{N-1}}\CO_ij\HES\notag\\
&=\ROm klij x_k\gamma_l +
(q-1)\left(1-q^{-1}\alpha\frac{1-q^{N-1}}{1+q^{N-2}}\right)x_i\gamma_j
\notag\\*&\quad +\frac{q^{N-2}(q^2-1)}{\alpha+1}
\left(\frac1{1-q^{N-1}}-\frac{q^{-1}\alpha}{1+q^{N-2}}\right)
\KO klijx_k\gamma_l\notag\\*&\quad +
\frac{q-1}{\alpha+1}\left(q^{-1}\alpha^2\frac{1-q^{N-1}}{1+q^{N-2}}
-q^{-1}\alpha(q-1) + \frac{1+q^{N-2}}{1-q^{N-1}}\right)
\CO^kl x_ix_jx_k\gamma_l\text.\label{eqbimg1}
\end{align}
The bimodule structure w.r.t.\ the new basis is therefore linear
if and only if the last coefficient is zero, i.e.
\begin{align*}
0&=\alpha^2-\frac{(q-1)(1+q^{N-2})}{1-q^{N-1}}
-q\frac{(1+q^{N-2})^2}{(1-q^{N-1})^2}\\
&=\left(\alpha+\frac{1+q^{N-2}}{1-q^{N-1}}\right)
\left(\alpha-q\frac{1+q^{N-2}}{1-q^{N-1}}\right)\text.
\end{align*}
The two solutions for $\alpha$ as stated in the
Proposition are now obvious;
by inserting them into equation \eqref{eqbimg1} one verifies
\eqref{eqbimg}.
This completes the proof for ${\Gamma=\GES+}$; the case
${\Gamma=\GES-}$ is treated in the same way.\qed

Since we have to buy the particularly
simple form of the bimodule structure in the $\gamma^{\pm}$
bases with more complicated equations describing the
differentiation map $\d$
\begin{equation*}
\d x_i=\gamma^-_i\mp q\frac{1+q^{N-2}}{1\pm q}x_i\CO^klx_k\gamma^-_l=
\gamma^+_i+q^{-N+2}\frac{1+q^{N-2}}{1\pm q}x_i\CO^klx_k\gamma^+_l\text,
\end{equation*}
it is wise to keep in mind both types of basis for
$\Gamma$ to use each one when appropriate, depending on the kind
of calculations to be done.

The $*$\nob-structure still takes a simple form in the new bases.
\begin{Cor}
Let ${\Gamma\in\{\GES+,\GES-\}}$. Then
${(\gamma^+_i)^*=q^{-1}\CO^ij\gamma^-_j}$;
${(\gamma^-_i)^*=q\CO^ij\gamma^+_j}$.
\end{Cor}
\textsc{Proof:} We prove the first equality for ${\Gamma=\GES+}$.
\begin{align*}
(\gamma^+_i)^*&=\d x^*_i - \frac{1+q^{N-2}}{1-q^{N-1}}\HES^*x^*_i=
\CO^ij\left(\d x_j + \frac{1+q^{N-2}}{1-q^{N-1}}\HES x_j\right)\\& =
\CO^ij\left(\d x_j + \frac{1+q^{N-2}}{1-q^{N-1}}x_j\HES +
(q^{-1}-1)\d x_j\right)= q^{-1}\CO^ij\gamma^-_j\text.
\end{align*}
The proofs of the remaining statements are analogous.\qed
\subsection{Proof of the classification statements}\label{sclpr}
The most important tool which will be employed in our proofs of
the classification theorem~\ref{thfcl} and Corollary~\ref{crffcl}
is to analyse morphisms of corepresentations of ${\cA=\Oq(n)}$
(co-) acting on the algebra ${X=\ESq{N-1}}$ and its tensor products.
In this way, covariance and the algebraic constraints are
exploited to find the general structure, with unknown coefficients,
for the bimodule structure of the desired differential calculi;
then the generating relations of the algebra and the definitions
for differential calculi are evaluated to determine the coefficients.

The same principle underlies the proofs given by Apel and Schm{\"u}dgen
in~\cite{as} (for Podle{\'s}' quantum spheres) and by the author
in~\cite{wk1} and~\cite{wk3} (for the Vaksman-Soibelman quantum spheres
and quantum projective spaces, resp.).

First we observe that since $q$ is not a root of unity,
the corepresentation theory for $\Oq(N)$ is quite similar to
that of the undeformed matrix group~\cite{hay} for which
the formalism of Young frames can be used~\cite{br}.

Denote by $V(k)$ the vector space of those polynomials of degree $k$
in the generators ${x_1,\dots,x_N}$ of $X$, and by $V_{\mathrm{r}}(k)$
the vector subspace of all those degree $k$ polynomials which
by virtue of the defining relations~\eqref{eqesr1} and~\eqref{eqesr2}
of~$\ESq{N-1}$ are equal to polynomials of lower degree.
Let $W(k)$ be the vector space complement of $V_{\mathrm{r}}(k)$
in $V(k)$. The coaction~$\DeltaR$ induces corepresentations
${\pi(k):W(k)\to W(k)\otimes\cA}$. In particular, $\pi(0)$ is
the trivial corepresentation ${1_X\mapsto1_X\otimes1_{\cA}}$,
denoted also by $\ynull$, and $\pi(1)$ is the fundamental
corepresentation, ${x_i\mapsto x_j\otimes u^j_i}$, which corresponds
to the Young frame $\ytt100010000$.

We can now use the Young frame formalism for the orthogonal
quantum group $\Oq(N)$ to calculate the higher
corepresentations~$\pi(k)$ and the tensor products
${\pi(k)\otimes\pi(1)}$ successively. Note that ${\pi(k+1)}$
is obtained by cancelling those summands of ${\pi(k)\otimes\pi(1)}$
which are annihilated by the commutation relation~\eqref{eqesr1}
or reduced to lower degree by the inhomogeneous relation~\eqref{eqesr2}.
\begin{align*}
\pi(0)\otimes\pi(1)&=\ytt100010000\text;\\
\pi(1)\otimes\pi(1)&=\ytt200011000+\ytt110020000+\ynull\text;&
\pi(2)&=\ytt200011000\text;\\
\pi(2)\otimes\pi(1)&=\ytt300011100+\ytt210021000+\ytt100010000\text;&
\pi(3)&=\ytt300011100\text;\\
\pi(3)\otimes\pi(1)&=\ytt400011110+\ytt310021100+\ytt200011000\text;&
\pi(4)&=\ytt400011110\text;\\
\pi(4)\otimes\pi(1)&=\ytt500011111+\ytt410021110+\ytt300011100\text;&
\pi(5)&=\ytt500011111\text.
\end{align*}
These equalities give the correct decompositions of $\pi(k)$ and
${\pi(k)\otimes\pi(1)}$ into irreducible corepresentations
for ${N\ge6}$; for smaller $N$ some of the summands do not
exist, or decompose. This is the
reason why we can guarantee completeness of our classification
only under this assumption.
For ${k\ge5}$, the corepresentations of both $\pi(k)$ and
${\pi(k)\otimes\pi(1)}$ contain only Young frames with 4 and more
columns.

In the language of corepresentations,
the classification constraint of Theorem~\ref{thfcl} says precisely
that $\Gamma$ is generated by ${\d V(1)}$ as a left module.
As a vector space with a corepresentation of $\Oq(N)$,
${\d V(1)}$ is the same as $V(1)$; therefore $\Gamma$ is the sum
of invariant vector spaces for all (irreducible components of)
corepresentations ${\pi(k)\otimes\pi(1)}$, ${k=0,1,2,\dots}$;
for the same reason, the vector subspace ${\d V(1)\cdot V(1)}$
of all ${\d x_i\cdot x_j}$
decomposes into the invariant vector spaces for the corepresentation
${\pi(1)\otimes\pi(1)}$.

The bimodule structure is described by a transformation that sends
one-forms from ${\d V(1)\cdot V(1)}$ to left-module expressions
in ${X\d V(1)}$.
The covariance condition requires that this transformation must
transform each invariant subspace for an irreducible corepresentation
of the quantum group to an isomorphic subspace. Thus, the transformation
is necessarily a linear combination of all morphisms from irreducible
invariant subspaces in ${\d V(1)\cdot V(1)\cong V(1)\otimes V(1)}$ to
isomorphic invariant subspaces in ${X\d V(1)\cong X\otimes V(1)}$.

Searching the sums of Young frames representing ${\pi(k)\otimes\pi(1)}$
for frames which are also contained in ${\pi(1)\otimes\pi(1)}$,
indicating isomorphic invariant subspaces, we find that
$\ytt200011000$ occurs in ${\pi(1)\otimes\pi(1)}$ and in
${\pi(3)\otimes\pi(1)}$ while $\ytt110020000$ and $\ynull$ occur
only in ${\pi(1)\otimes\pi(1)}$ itself. Therefore one has
four relevant morphisms and gets the ansatz
\begin{align*}
\d x_i\cdot x_j
&=\alpha_1(P_+)^{kl}_{ij}x_k\d x_l
+\alpha_2(P_-)^{kl}_{ij}x_k\d x_l\\*&\quad\qquad
+\alpha_3(P_0)^{kl}_{ij}x_k\d x_l
+\alpha_4(P_+)^{st}_{ij}\CO^klx_sx_tx_k\d x_l
\end{align*}
for the bimodule structure of any free covariant first order
differential calculus on $\ESq{N-1}$.
Here, $P_+$ is the projector onto the subspace corresponding
to $\ytt200011000$, $P_-$ the projector onto the subspace
corresponding to $\ytt110020000$, and $P_0$ the projector
onto the subspace for $\ynull$. Since the projectors are less
convenient for later calculation, we use the fact that
$\ROm{}{}{}{}$, the identity $\IO {}{}{}{}$, and $\KO {}{}{}{}$
are three independent linear combinations of $P_+$, $P_-$ and
$P_0$, spanning therefore the same vector space of transformations,
and that the multiplication $x_sx_t$ in the last summand in fact
does itself the job of $P_+$ because it involves a symmetrisation
via relations~\eqref{eqesr1} and~\eqref{eqesr2}. We can therefore
rewrite the ansatz as
\begin{equation}\label{eqfabs}
\d x_i\cdot x_j=a_1\ROm klij x_k\d x_l+a_2x_i\d x_j
+a_3\KO klij x_k\d x_l+a_4\CO^kl x_ix_jx_k\d x_l
\end{equation}
with unknown coefficients $a_1$, $a_2$, $a_3$ and~$a_4$.

The changed constraint in Corollary~\ref{crffcl} requiring the
calculus to contain the left-module relations generated by $\HES=0$
means that the subspace for $\ynull$ is cancelled from
${\pi(1)\otimes\pi(1)}$ and that for $\ytt200011000$ from
${\pi(3)\otimes\pi(1)}$ in the decomposition of~${X\d V(1)}$.
This reduces the number of possible morphisms to~2, and with
essentially the same arguments as before the ansatz
\begin{equation}\label{eqffabs}
\d x_i\cdot x_j=a_1\ROm klij x_k\d x_l+a_2x_i\d x_j
\end{equation}
is obtained.

We return now to equation~\eqref{eqfabs}. To determine the
coefficients we use first the conditions
\begin{align}
\CO^ij(\d x_i\cdot x_j+x_i\d x_j)&=0\text,\label{eqca1}\\
\ROm klij(\d x_k\cdot x_l+x_k\d x_l)-q(\d x_i\cdot x_j+x_i\d x_j)&=0\text,
\label{eqca2}\\
\CO^kl\d x_i\cdot x_kx_l&=\d x_i\text,\label{eqca3}\\
\RO stkl\d x_i\cdot x_sx_t-q\d x_i\cdot x_kx_l+
\frac{q-q^{-1}}{1+q^{N-2}}\CO_kl\d x_i&=0\text,\label{eqca4}
\end{align}
the first two of which arise by differentiating the defining
relations of $\ESq{N-1}$ while the last two ones result from
the bimodule requirement and the defining relations of $\ESq{N-1}$.
Since covariance and freeness are guaranteed by the construction
leading to~\eqref{eqfabs}, these equations which completely encode
the compatibility of bimodule structure with differentiation are
everything which remains to be satisfied.

Each of the equations~\eqref{eqca1}--\eqref{eqca4} is evaluated
by transforming it via~\eqref{eqfabs} to left-module form%
\footnote{Part of these calculations has been carried out with the
aid of a computer algebraic program.}
and comparing coefficients
for elements of ${X\d V(1)}$.%
\footnote{The coefficient comparisons involved
are admissible only if there are enough algebraically independent
elements in ${X\d V(1)}$ but this requirement is met for ${N\ge6}$.}

From~\eqref{eqca1} and~\eqref{eqca2} we derive
\begin{align*}
a_2&=qa_1-1;&a_4&=-1-q^{N-1}a_1-a_2
-\frac{q^2(1+q^{N-2})(1-q^{-N})}{q^2-1}a_3
\end{align*}
and eliminate $a_2$ and $a_4$. After doing so, we obtain from
\eqref{eqca3} and~\eqref{eqca4} among others the conditions
\begin{gather*}
{a_1}^2-1=0\text;\\-1-a_3-q^{N-3}\frac{q^2-1}{1+q^{N-2}}a_1
+q^{-1}\frac{1+q^N}{1+q^{N-2}}a_1a_3+\frac{1+q^N}{1+q^{N-2}}{a_1}^2
=0\text.
\end{gather*}
By computing for each of the two possible values for $a_1$ the
other coefficients, we find that two sets of coefficients are
still possible, namely
\begin{align*}
a_1&=1\text,&a_2&=q-1\text,&a_3&=q^{N-2}\frac{q^2-1}{1-q^{N-1}}\text,&
a_4&=\frac{(1-q)(1+q^{N-2})}{1-q^{N-1}}\text;\\
a_1&=-1\text,&a_2&=-q-1\text,&a_3&=q^{N-2}\frac{q^2-1}{1-q^{N-1}}\text,&
a_4&=\frac{(1+q)(1+q^{N-2})}{1+q^{N-1}}\text.
\end{align*}
Inserting these into~\eqref{eqfabs} yields equation~\eqref{eqbimx}
from the theorem for the $\GES+$ and $\GES-$ cases, resp.
To finish the proof of the theorem, one checks that
\eqref{eqca1}--\eqref{eqca4} are satisfied
by both sets of coefficients. For this last step, the assumption
${N\ge6}$ is not needed.

To prove Corollary~\ref{crffcl},
conditions~\eqref{eqca2}--\eqref{eqca4}
have to be exploited in the same manner but with the simpler
assumption~\eqref{eqffabs}. Equation~\eqref{eqca1} is trivially
satisfied in this case. The calculations are essentially the same
as above but simpler. Equation~\eqref{eqca2} leads again to
${a_2=qa_1-1}$; from~\eqref{eqca4} we get then, among other
conditions,
\begin{align*}
{a_1}^2&=1&&\text{and}&{a_1}^2-(q+q^{-1})a_1+1&=0
\end{align*}
which contradict each other except for ${q=\pm1}$. This completes
the proof.
\section{Higher order differential calculus and symmetry}
\subsection{Basic definitions}
If the geometric structure of quantum spaces is to be investigated,
first order differential calculus is an insufficient tool since
even simple concepts of differential geometry require at least
second order differential forms to be formulated. That's why we
shall turn our interest now to higher order differential calculus.

Given a covariant first order differential calculus $\Gammad$
over a quantum homogeneous space $X$, consider
a pair ${(\Gw,\d)}$ consisting of a graded algebra
$\Gw$ and a linear mapping ${\d:\Gw\to\Gw}$
such that the degree $0$ and $1$ components of $\Gw$
are isomorphic to~$X$ and~$\Gamma$, resp. Let the multiplication of
$\Gw$ be denoted by $\wedge$, with the convention that
the $\wedge$ sign may be omitted if one of the factors involved
is of degree zero (i.e.\ an element of~$X$). Suppose further that
the following statements hold for ${(\Gw,\d)}$:
The mapping $\d$ increases the degree by $1$;
$\d$ extends the differential of the first order calculus and
fulfils the graded Leibniz rule,
${\d(\vartheta_1\wedge\vartheta_2)=}{\d\vartheta_1\wedge\vartheta_2
+(-1)^d\vartheta_1\wedge\d\vartheta_2}$,
for ${\vartheta_1,\vartheta_2\in\Gw}$
where $d$ is the degree of $\vartheta_1$; and ${\d\d=0}$.
Finally, assume that the covariance map $\PhiR$
from the first order calculus can be extended to a map
${\PhiR^{\wedge}:\Gw\to\Gw\otimes{\cA}}$
making $\Gw$ into a covariant $X$-bimodule.
If all these properties are fulfilled, the pair ${(\Gw,\d)}$
is called a \textsl{covariant higher order differential
calculus over~$X$ which extends~$\Gammad$.}

First, one observes that for any first order differential calculus
$\Gamma$ over a quantum homogeneous space, there exists a
universal higher order differential calculus $\Gw[u]$ extending
$\Gamma$ from which any other calculus with the above properties
can be obtained by factorisation. Relations of this calculus
can be successively computed started with those of the first
order differential calculus and applying the definition of
higher order differential calculus but in general the universal
calculus is difficult to handle and displays properties which
make not much sense for geometric investigations, such as unlimited
order of differential forms on deformations of finite-dimensional
manifolds. Therefore it is mostly of merely algebraic interest.

A second approach to higher order differential calculus is based
on an antisymmetrisation procedure, cf.~\cite{wor},~\cite{sch}.
For this, a ``braiding'' map
${\sigma:\Gamma\otimes_{X}\Gamma\to\Gamma\otimes_{X}\Gamma}$
is required which must be an algebra homomorphism satisfying
the braid (or quantum Yang-Baxter) equation
\begin{equation*}
(\sigma\otimes\id)(\id\otimes\sigma)(\sigma\otimes\id)=
(\id\otimes\sigma)(\sigma\otimes\id)(\id\otimes\sigma)\text.
\end{equation*}
The algebra $\Gw$ is then obtained by factorising $\Gof{}$
over an ideal which is generated by the kernel of
${(\id_{\Gamma\otimes\Gamma}-\sigma)}$.

In the literature several versions for the construction of
the ideal can be found. All of them agree in their second order
component but differences occur in higher orders, leading possibly to
different external algebras. For a detailed comparison in the
case of the quantum group $\mathrm{SL}_q(N)$ see~\cite{sch}.

An advantage of the antisymmetrisation approach lies in the
symmetry information encoded in $\sigma$ which merits an
interest of its own in further exploring of the differential
structure.

\subsection{Second order relations in the universal higher order
differential calculus extending \mbm{$\GES+$}}
For the remaining part of this article,
we shall base our considerations on the first order differential
calculus ${\Gamma=\GES+}$ over ${X=\ESq{N-1}}$.

The first order differential calculus determines a unique
universal higher order differential calculus $\Gw[u]$ over
$\ESq{N-1}$ which extends $\GES+$.
In the following we want to give relations which have to hold
in this higher order differential calculus. In particular,
we want to give a full account of the relations of second order.

Since iterated differentiation is involved in the calculations,
we prefer to use the ${\d x_i}$ basis for $\Gamma$.
\begin{Prop}\label{prhu}
All left-module relations in the module ${\Gamma\wedge\Gamma}$
of two-forms in the universal higher order differential calculus
on $\ESq{N-1}$ extending $\GES+$ are generated by the set of
relations
\begin{align}
0&=\ROm klij \d x_k\wedge\d x_l+ q\d x_i\wedge\d x_j
- q\frac{(1-q)(1+q^{N-2})}{1-q^{N-1}}x_i\HES\wedge\d x_j\notag\\&\quad
-\! \frac{(1-q)(1+q^{N-2})}{1-q^{N-1}}\ROm klij x_k\HES\wedge\d x_l\!
-\! \frac{(1+q^2)(1+q^{N-2})^2}{(1-q^{N-1})^2}x_ix_j\CO^klx_k\HES\wedge\d x_l
\notag\\&\quad+ q^{N-2}\frac{(1+q)(1+q^2)(1+q^{N-2})}{(1-q^{N-1})^2}
\KO klij x_k\HES\wedge\d x_l\text,\label{eqiw}\\
0&=\d\HES+2q\frac{1+q^{N-2}}{(1-q)(1-q^{N-1})}\CO^kl x_k\HES\wedge\d x_l
\text.\label{eqiiw}
\end{align}
\end{Prop}
\textsc{Proof:}
We start by differentiating equations \eqref{eqbimx} for $\GES+$
using the graded Leibniz rule. By transforming all involved
summands to left-module expressions, we obtain a relation
which necessarily holds in $\Gw[u]$.
\begin{align*}
0&=   \ROm klij \d x_k\wedge\d x_l+ q   \d x_i\wedge\d x_j
- q^{N-1} \frac{(1-q)^2(1+q)(1+q^{N-2})}{(1-q^{N-1})^2}x_i\HES\wedge\d x_j
\\&\quad+ q^2 \frac{(1-q)(1+q^{N-2})}{1-q^{N-1}}x_ix_j\d\HES
+ q \frac{(1-q)(1+q^{N-2})}{1-q^{N-1}}\ROm stjk\CO^kl x_ix_s
\d x_t\wedge\d x_l\\&\quad
- q^{N-2} \frac{(1-q)^2(1+q)(1+q^{N-2})}{(1-q^{N-1})^2}
\ROm klij x_k\HES\wedge\d x_l\\&\quad+ \frac{(1-q)(1+q^{N-2})}{1-q^{N-1}}
\ROm tukl\ROm skij\CO^lv x_sx_t\d x_u\wedge\d x_v\\&\quad
- q^{N-2} \frac{(1+q)(1-q)}{1-q^{N-1}}\CO_ij \d\HES\!
+\! \frac{(1-q)(1-q^3)(1+q^{N-2})^2}{q(1-q^{N-1})^2}
\CO^kl x_ix_jx_k\HES\wedge\d x_l\\*&\quad
- q^{N-3} \frac{(1-q)^2(1+q)(1+q^{N-2})}{(1-q^{N-1})^2}
\KO klij x_k\HES\wedge\d x_l
\end{align*}

This equation is resolved for ${\ROm klij \d x_k\wedge\d x_l}$; by
substituting the resulting expression on the right-hand side of the
identity
\begin{equation*}
\ROm stjk\CO^kl x_ix_s\d x_t\wedge\d x_l
=\bigl(\ROm tlsj\CO^ks+(q-q^{-1})(\IO tlsj-\KO tlsj)\CO^ks\bigr)
x_ix_s\d x_t\wedge\d x_l\text,
\end{equation*}
a new equation is obtained which contains again
${\ROm stjk\CO^kl x_ix_s\d x_t\wedge\d x_l}$ on its right-hand side.
Resolving for this term, one finds
\begin{align*}
\ROm stjk\CO^kl x_ix_s\d x_t\wedge\d x_l&=
-\frac{1-q^{N-2}-q^{N-1}+q^N}{1-q^{N-1}}x_i\HES\wedge\d x_j
+ (1-q) x_ix_j\d\HES\\*&\quad\quad
- \frac{(1-q)(1+q^{N-2})}{1-q^{N-1}}\CO^kl x_ix_jx_k\HES\wedge\d x_l\text,
\end{align*}
which simplifies the original relation to
\begin{align*}
0&=   \ROm klij \d x_k\wedge\d x_l+ q   \d x_i\wedge\d x_j
+ q \frac{(1-q)(1+q^{N-2})}{1-q^{N-1}}x_i\HES\wedge\d x_j\\*&\quad
- q^{-1} \frac{(1-q)(1-q+q^2)(1+q^{N-2})}{1-q^{N-1}}x_ix_j\d\HES
\\*&\quad+ \frac{(1-q)(1+q^{N-2})}{1-q^{N-1}}\ROm klij x_k\HES\wedge\d x_l
\\*&\quad+ q^{N-3} \frac{(1-q)(1+q^3)}{1-q^{N-1}}\CO_ij \d\HES
- \frac{(1-q)^2(1+q^{N-2})^2}{(1-q^{N-1})^2}\CO^kl x_ix_jx_k\HES\wedge\d x_l
\\*&\quad+ q^{N-2} \frac{(1-q)^2(1+q)(1+q^{N-2})}{(1-q^{N-1})^2}
\KO klij x_k\HES\wedge\d x_l\text.
\end{align*}
\newcommand{\eR}[2]{\mathbf{[R_{\mathnormal{#1#2}}]}}%
We shall denote the right-hand side of this equation,
with the $\wedge$ signs replaced by $\otimes_X$ (and thus
$\d\HES$ by ${\d^\otimes\HES:=\CO^ij\d x_i\otimes_X\d x_j}$),
by $\eR ij$.

By factorising the free left $\ESq{N-1}$\nob-module
${\Gamma\otimes_X\Gamma}$ generated by ${\d x_i\otimes\d x_j}$
over the relations
$\eR ij=0$, a left $\ESq{N-1}$\nob-module $\Gamma^2_1$ is obtained.
This is not a bimodule, and therefore still not
${\Gamma\wedge\Gamma}$, because if this were the case,
it would follow from equation \eqref{eqbimx} and $\eR ij=0$ that
\newcommand{\eRR}[3]{\mathbf{[R'_{\mathnormal{#1#2#3}}]}}
\begin{align*}
0&=\eR ij x_k=\eRR ijk\\*&=  A_1 x_ix_jx_k\d^\otimes\HES
+ A_2 \CO_ij x_k\d^\otimes\HES+ A_3 \CO_jk x_i\d^\otimes\HES
+ A_4 \CO_tk\ROm stij x_s\d^\otimes\HES\\*&\quad
+ A_5 \CO^mn x_ix_jx_kx_m\HES\otimes_X\d x_n
+ A_6 \CO_ij \CO^mn x_kx_m\HES\otimes_X\d x_n\\*&\quad
+ A_7 \CO_jk \CO^mn x_ix_m\HES\otimes_X\d x_n
+ A_8 \CO_tk\ROm stij\CO^mn x_sx_m\HES\otimes_X\d x_n\text,
\end{align*}
where
\begin{align*}
A_5 &= -2q^{-3} \frac{\begin{array}{l@{}l}(1+q^{N-2})^2 &
(1-2q+2q^2-2q^3+2q^4-2q^5+q^6-q^{N-1}\\~&~~
+2q^N-2q^{N+1}+q^{N+2}+2q^{N+4}-3q^{N+5}+q^{N+6})\end{array}
}{(1-q^{N-1})^3}\text,\\A_6 &= 2q^{N-5} \frac{\begin{array}{l@{}l}
(1+q^{N-2}) &(1-2q+2q^3-q^5-q^6+q^7-q^{N-1}+2q^N\\~&~~
-3q^{N+2}+q^{N+3}+3q^{N+4}-q^{N+5}-2q^{N+6}+q^{N+7})\end{array}
} {(1-q^{N-1})^3}\text,\\A_7 &= qA_8 = 2q^{N-3}
\frac{(1-q)^4(1+q)^2(1+q^{N-2})}{(1-q^{N-1})^2}\text,\\
\frac{A_5}{A_1} &= \frac{A_6}{A_2}=\frac{A_7}{A_3}=\frac{A_8}{A_4}
=T:=2q\frac{1+q^{N-2}}{(1-q)(1-q^{N-1})}\text.
\end{align*}
To make $\Gamma^2_1$ into a bimodule,
it has to be factorised again over the relation
\begin{equation}\eRR ijk=0\text.\label{equhr2}\end{equation}
We observe that there is an ${S\ne0}$ such that
\begin{equation*}
\CO^ij
(A_1 x_ix_jx_k + A_2\CO_ij x_k + A_3\CO_jk x_i + A_4\CO_tk\ROm stij x_s)
= S x_k\text.
\end{equation*}
Because of this fact and the equality
\begin{align*}
&\eRR ijk=\\*
&(A_1 x_ix_jx_k + A_2\CO_ij x_k + A_3\CO_jk x_i + A_4\CO_tk\ROm stij x_s)
(\d^\otimes\HES+T\CO^mn x_m\HES\otimes\d x_n)\text,
\end{align*}
the relation \eqref{equhr2} implies
\begin{equation*}
0=\CO^rk\CO^ij x_r\eRR ijk
\end{equation*}
and thus
\begin{equation}
\d^\otimes\HES+T\CO^mn x_m\HES\otimes_X\d x_n=0\text.
\label{equhr3}
\end{equation}
Let $\Gamma^2_2$ be the left module obtained by factorising
$\Gamma^2_1$ over the relation \eqref{equhr3}.

Since in $\Gamma^2_2$ the equalities
\begin{align*}
\eR ij x_k&=0&&\text{and}&
(\d^\otimes\HES+T\CO^mn x_m\HES\otimes_X\d x_n)x_k&=0
\end{align*}
are satisfied, $\Gamma^2_2$ is in fact a $\ESq{N-1}$\nob-bimodule
and therefore isomorphic to ${\Gamma\wedge\Gamma}$.\qed
\subsection{Braiding symmetry and antisymmetrisation}
From the simple bimodule structure of $\GES+$ w.r.t.\ the basis
$\{\gamma_i\}$ where either ${\gamma_i=\gamma^+_i}$ for ${i=1,\dots,N}$,
or ${\gamma_i=\gamma^-_i}$ for ${i=1,\dots,N}$, one easily guesses
the ansatz for a covariant ``braiding'' homomorphism which underlies
the following theorem.

\begin{Thm}\label{thsy}
Let ${X=\ESq{N-1}}$, ${\Gamma=\GES+}$, and either
${\gamma_i=\gamma^+_i}$ for ${i=1,\dots,N}$ or
${\gamma_i=\gamma^-_i}$ for ${i=1,\dots,N}$.
For any real number $\alpha$, the equation
\begin{equation*}
{\sigma(\gamma_i\otimes_X\gamma_j)}=
{\alpha\ROm klij\gamma_k\otimes_X\gamma_l}
\end{equation*}
continues to a well-defined covariant homomorphism
${\sigma:\Gamma\otimes_X\Gamma\to\Gamma\otimes_X\Gamma}$
which satisfies the braid relation.
\end{Thm}
\textsc{Proof:}
Since all tensor products of differential modules in the
following have to be taken over the quantum space $\ESq{N-1}$,
we shall often write $\otimes$ instead
of $\otimes_X$ or $\otimes_{\ESq{N-1}}$ in the following.

First, let ${\syL:\Gamma\otimes\Gamma\to\Gamma\otimes\Gamma}$
be the left module homomorphism defined by
${\syL(\gamma_i\otimes\gamma_j)=}{\alpha\ROm klij\gamma_k\otimes\gamma_l}$.
Because of
\begin{align*}
&\syL(\gamma_i\otimes\gamma_j\cdot x_k)
-\syL(\gamma_i\otimes\gamma_j)\cdot x_k\\
&\qquad=\syL(\ROm stiu\ROm uvjk x_s\gamma_t\otimes\gamma_v)
-\alpha\ROm uvij\gamma_u\otimes\gamma_v\cdot x_k\\
&\qquad=\alpha(\ROm wptv\ROm stiu\ROm uvjk
-\ROm swut\ROm tpvk\ROm uvij) x_s\gamma_w\otimes\gamma_p =0\text,
\end{align*}
$\syL$ is also a right module homomorphism, i.e.\ $\sigma$ is a
well-defined bimodule homomorphism.

Second, with the abbreviations
${\sigma_1:=\sigma\otimes_X\id_\Gamma}$ and
${\sigma_2:=\id_\Gamma\otimes_X\sigma}$, we obtain
\begin{align*}
&(\sigma_1\sigma_2\sigma_1-\sigma_2\sigma_1\sigma_2)
(\gamma_i\otimes\gamma_j\otimes\gamma_k)\\
&\quad=\alpha^3(\ROm swut\ROm tpvk\ROm uvij
-\ROm wptv\ROm stiu\ROm uvjk)
\gamma_s\otimes\gamma_w\otimes\gamma_p=0\text,
\end{align*}
i.e.\ the braid relation is fulfilled.

Finally, for the mapping
${\PhiR\otimes\PhiR:\Gamma\otimes_X\Gamma\to\Gamma\otimes_X\Gamma
\otimes\Oq(N)}$ (where the multiplication of $\Oq(N)$
makes the second tensor factors of both $\PhiR$ images into
one element of the quantum group) one has
\begin{align*}
&(\PhiR\otimes\PhiR)\sigma(\gamma_i\otimes\gamma_j)=
\alpha\ROm klij\gamma_s\otimes\gamma_t\otimes u^s_ku^t_l\\
&\quad=\alpha\ROm stkl\gamma_s\otimes\gamma_t\otimes u^k_iu^l_j
=(\sigma\otimes\id_\Gamma)(\PhiR\otimes\PhiR)
(\gamma_i\otimes\gamma_j)\text,
\end{align*}
which guarantees the covariance of $\sigma$.\qed

\textsl{Remark:} Obviously, $\sigma$ is invertible, with
the inverse $\sigma^{-1}$ defined by
\begin{equation*}
{\sigma^{-1}(\gamma_i\otimes\gamma_j)}=
{\alpha^{-1}\RO klij\gamma_k\otimes\gamma_l}\text.
\end{equation*}
\begin{Cor}
The homomorphism $\sigma$ constructed in Theorem~\ref{thsy}
is independent on the choice ${\gamma_i=\gamma^+_i}$ or
${\gamma_i=\gamma^-_i}$.
\end{Cor}
\textsc{Proof:}
Define $\sigma^-$ and $\sigma_+$ via
\begin{equation*}
{\sigma^-(\gamma^-_i\otimes\gamma^-_j)}=
{\alpha\RO klij\gamma^-_k\otimes\gamma^-_l}\text,\quad
{\sigma^+(\gamma^+_i\otimes\gamma^+_j)}=
{\alpha\RO klij\gamma^+_k\otimes\gamma^+_l}\text.
\end{equation*}
From Proposition~\ref{prga} and the bimodule structure of $\GES+$
it follows that
\begin{align*}
\gamma^+_i\otimes\gamma^+_j&=\gamma^-_i\otimes\gamma^-_j
-(1+q^{N-2})\CO^klx_ix_k\gamma^-_k\otimes\gamma^-_l\\*&\quad\quad
-(1+q^{N-2})\CO^vt\ROm tsuv\ROm kuij x_kx_l\gamma^-_s\otimes\gamma^-_t
\\*&\quad\quad+q^2(1+q^{N-2})^2\CO^kl\CO^stx_ix_jx_kx_s\gamma^-_t
\otimes\gamma^-_l\\*&\quad\quad-(q-q^{-1})(1+q^{N-2})\CO^ul\ROm skju
x_ix_s\gamma^-_k\otimes\gamma^-_l\\*&\quad\quad
-(q^2-1)(1+q^{N-2})\CO^klx_ix_j\gamma^-_k\otimes\gamma^-_l\text.
\end{align*}
From this one can calculate
${\alpha^{-1}\sigma^-(\gamma^+_i\otimes\gamma^+_j)}$ as well as
${\ROm klij\gamma^+_k\otimes\gamma^+_l}$. One finds that both
equal the expression
\begin{align*}
&\ROm klij\gamma^-_k\otimes\gamma^-_l
-(1+q^{N-2})\CO^su\ROm kluj x_ix_s\gamma^-_l\otimes\gamma^-_l\\*
&-(1+q^{N-2})\CO^kl\ROm stij x_sx_k\gamma^-_l\otimes\gamma^-_t
+q(1+q^{N-2})^2\CO^kl\CO^stx_ix_jx_kx_s\gamma^-_t\otimes\gamma^-_l\\*
&-(q-q^{-1})(1+q^{N-2})\CO^klx_ix_k\gamma^-_l\otimes\gamma^-_j
\end{align*}
which means that $\sigma^+$ and $\sigma^-$ are identical.\qed

The following corollary specifies the way in which
$\sigma$ can be used to define a higher order differential
calculus extending $\GES+$ via antisymmetrisation.
\def\eI{\mathbf{[I]}}\def\eII{\mathbf{[II]}}%
Note that the expressions $\eI$ and $\eII$ correspond to
the right-hand sides of the relations \eqref{eqiw}, \eqref{eqiiw}
from Proposition~\ref{prhu}
\begin{Cor}
Let $X$, $\Gamma$, $\gamma_i$, $\alpha$ and $\sigma$ be as in
Theorem~\ref{thsy}. If and only if ${\alpha=q}$, the homomorphism
${(\id_{\Gamma\otimes\Gamma}-\sigma)}$ annihilates the two expressions
${\eI,\eII\in\Gamma\otimes\Gamma}$ defined by
\begin{align*}
\eI&:=\ROm klij \d x_k\otimes\d x_l+ q\d x_i\otimes\d x_j\\*&\quad\quad
- q\frac{(1-q)(1+q^{N-2})}{1-q^{N-1}}x_i\CO^mn x_m\d x_n\otimes\d x_j
\\*&\quad\quad- \frac{(1-q)(1+q^{N-2})}{1-q^{N-1}}\ROm klij x_k
\CO^mn x_m\d x_n\otimes\d x_l\\*&\quad\quad
- \frac{(1+q^2)(1+q^{N-2})^2}{(1-q^{N-1})^2}x_ix_j
\CO^klx_k\CO^mn x_m\d x_n\otimes\d x_l\\*&\quad\quad
+ q^{N-2}\frac{(1+q)(1+q^2)(1+q^{N-2})}{(1-q^{N-1})^2}
\KO klij x_k\CO^mn x_m\d x_n\otimes\d x_l\text,\\
\eII&:=\CO^mn\d x_m\otimes\d x_n+2q\frac{1+q^{N-2}}{(1-q)(1-q^{N-1})}\CO^kl
\CO^mn x_kx_m\d x_n\otimes\d x_l\text.
\end{align*}
\end{Cor}
\textsc{Proof:}
That the parameter $\alpha$ must be chosen as $q$ in order
for $\eI$ and $\eII$ to be annihilated by ${(\id-\sigma)}$
is obvious. Let therefore
\begin{equation*}
\sigma(\gamma_i\otimes_X\gamma_j)=
q\ROm klij\gamma_k\otimes_X\gamma_l\text.
\end{equation*}
By calculation, we obtain successively
\begin{align*}
&\d x_i\otimes_X\d x_j=\gamma_i\otimes_X\gamma_j- q\frac{1+q^{N-2}}{1+q}
\CO^mn x_ix_m\gamma_n\otimes_X\gamma_j\\*&\quad\quad
+ q^2\frac{(1-q)(1+q^{N-2})}{1+q}\CO^mnx_ix_j\gamma_m\otimes_X\gamma_n
\\&\quad\quad+ q\frac{(1-q)(1+q^{N-2})}{1+q}
\ROm kljm\CO^mnx_ix_k\gamma_l\otimes_X\gamma_n\\&\quad\quad
- q\frac{1+q^{N-2}}{1+q}\ROm kltm\ROm stij\CO^mnx_sx_k\gamma_l\otimes_X
\gamma_n\\*&\quad\quad+ q^4\frac{(1+q^{N-2})^2}{(1+q)^2}
\CO^mn\CO^stx_ix_jx_sx_m\gamma_n\otimes_X\gamma_t\text,\\
&\eI=q\gamma_i\otimes_X\gamma_j+\ROm klij\gamma_k\otimes_X\gamma_l
- q\frac{1+q^{N-2}}{1+q}\CO^mn x_ix_m\gamma_n\otimes_X\gamma_j\\*&\quad\quad
+ q^{-1}\frac{(1-q)(1+q^2)(1+q^{N-2})}{1-q^{N-1}}
\CO^mnx_ix_j\gamma_m\otimes_X\gamma_n\\&\quad\quad
- q^2\frac{1+q^{N-2}}{1+q}\ROm kljm\CO^mnx_ix_k\gamma_l\otimes_X\gamma_n
\\&\quad\quad- \frac{1+q^{N-2}}{1+q}\ROm klij\CO^mnx_kx_m\gamma_n\otimes_X
\gamma_l\\&\quad\quad- q\frac{1+q^{N-2}}{1+q}
\ROm kltm\ROm stij\CO^mnx_sx_k\gamma_l\otimes_X\gamma_n\\&\quad\quad
- q^{N-3}\frac{(1-q)(1+q^3+q^N)}{1-q^{N-1}}\KO klij\gamma_k\otimes_X\gamma_l
\\&\quad\quad- \frac{(1-q)(1+q^2(1+q^{N-2})^2}{(1+q)(1-q^{N-1})}
\CO^mn\CO^stx_ix_jx_sx_m\gamma_n\otimes_X\gamma_t\\*&\quad\quad
+ q^{N-2}\frac{(1-q+2q^2-q^3-q^{N+1})(1+q^{N-2})}{1-q^{N-1}}
\KO klij\CO^mn x_kx_m\gamma_n\otimes_X\gamma_l\text,\\[1ex]
&\CO^ij\d x_i\otimes_X\d x_j=\frac{2-q^{N-1}+q^N}{1+q}
\CO^ij\gamma_i\otimes_X\gamma_j\\*&\quad\quad
- q\frac{(1+q^{N-2})(1+q^N)}{(1+q)^2}
\CO^ij\CO^mnx_ix_m\gamma_n\otimes_X\gamma_j\text,\\
&\eII=- q^{N-1}\CO^ij\gamma_i\otimes_X\gamma_j+ q\frac{(1+q^{N-2})(1-q^N)}
{1-q^2}\CO^ij\CO^mnx_ix_m\gamma_n\otimes_X\gamma_j\text,
\end{align*}
and finally
\begin{align*}
(\id_{\Gamma\otimes\Gamma}-\sigma)\eI&=0\text,&
(\id_{\Gamma\otimes\Gamma}-\sigma)\eII&=0\text.
\end{align*}
\qed

It follows that a higher order differential calculus
extending $\GES+$ is obtained
by factorising $\Gof{}$ over the ideal generated by the kernel of
${(\id_{\Gamma\otimes\Gamma}-\sigma)}$ if and only if ${\alpha=q}$
which will be assumed from now on.

Let $\Gw[\sigma]$ be this higher order differential calculus.
Although $\Gw[\sigma]$ is smaller than $\Gw[u]$, it still
contains nonvanishing differential forms of arbitrarily
high order, e.g.\ ${\d\HES\wedge\d\HES\wedge\dots\wedge\d\HES}$.

In fact, ${(\id-\sigma)}$ cancels the component from
${\Gamma\otimes\Gamma}$ belonging to projector $P_+$ but leaves
the $P_0$ component intact. To complete antisymmetrisation,
the latter has to be zeroed by an additional factorisation.
Obviously, factorisation does not destroy the property of
being a higher order differential calculus.
\begin{Cor}
Let the higher order differential calculus $\Gw[\sigma0]$
be given by factorisation of $\Gw[\sigma]$
over the additional relation~${\d\HES=0}$. Then
$\Gw[\sigma0]$ is generated as a left $X$\nob-module by
the set of all differential forms
${\gamma_{i_1}\wedge\dots\wedge\gamma_{i_s}}$ where
${0\le s\le N}$ and ${1\le i_1<\dots<i_s\le N}$.
\end{Cor}
\textsc{Proof:}
By virtue of the bimodule structure of $\Gamma$,
any differential form in~$\Gw[\sigma0]$ can be
written as a sum of products
${x\gamma_{i_1}\wedge\gamma_{i_2}\wedge\dots\wedge\gamma_{i_s}}$
with ${x\in X}$, ${s\ge0}$ and ${i_1,i_2,\dots,i_s\in\{1,\dots,N\}}$.
Moreover, from ${\d\HES=0}$ it follows that
\begin{align*}
\CO^ij\gamma_i\wedge\gamma_j&=0\text,
&\CO^ij\CO^klx_ix_k\gamma_l\wedge\gamma_j&=0\text.
\end{align*}
Using
\begin{equation*}
(\id_{\Gamma\otimes\Gamma}-\sigma)
(\ROm klij\gamma_k\otimes\gamma_l+q\gamma_i\otimes\gamma_j)=
-q^{N-1}(q^2-1)\KO klij\gamma_k\otimes\gamma_l
\end{equation*}
we obtain therefore
\begin{equation}
\ROm klij\gamma_k\wedge\gamma_l+q\gamma_i\wedge\gamma_j=0\text.
\label{epw1}
\end{equation}
This commutation relation for the $\gamma_i$s allows to
replace any product ${\gamma_i\wedge\gamma_j}$, ${i>j}$ with a sum of
products ${\gamma_k\wedge\gamma_l}$ with ${k<l}$ and ${k+l=i+j}$.
Furthermore, the relation implies that the product
${\gamma_i\wedge\gamma_i}$ is zero if ${2i\ne N+1}$
while for ${2i=N+1}$ it can be replaced by a sum of
products ${\gamma_k\wedge\gamma_l}$ with ${k<l}$ and again ${k+l=i+j}$.

Assume now we are given the product
${\gamma_{i_1}\wedge\dots\wedge\gamma_{i_r}\wedge
\gamma_{i_{r+1}}\wedge\dots\wedge\gamma_{i_s}}$
where ${i_r\ge i_{r+1}}$ holds for some $r$.
Applying \eqref{epw1} as described before at positions $r$ and ${r+1}$
we transform the given product either to zero or to a sum of new products
${\gamma_{i'_1}\wedge\dots\wedge\gamma_{i'_r}\wedge
\gamma_{i'_{r+1}}\wedge\dots\wedge\gamma_{i'_s}}$
where ${i_k=i'_k}$ for all ${k\ne r,r+1}$ such that
${i'_1+i'_2+\dots+i'_k=i_1+i_2+\dots+i_k}$ whenever ${k\ne r}$
while ${i'_1+i'_2+\dots+i'_r\le i_1+i_2+\dots+i_r-1}$.
Since in each resulting product in each step one of the index sums
${i'_1+i'_2+\dots+i'_k}$ is by at least one smaller than the
corresponding sum in its predecessor, and since none of the
index sums can become negative, the procedure terminates after a
finite number of steps, yielding only products with ${i_r<i_{r+1}}$
for all $r$. This completes the proof.\qed
\section*{Acknowledgement}
I want to thank Dr.\ Istv{\'a}n Heckenberger for
many helpful discussions.

Part of the work included in this paper was supported by
Deutsche Forschungsgemeinschaft within the programme of
Graduiertenkolleg Quantenfeldtheorie, Leipzig.

\newcommand{\mbi}[2][.]{\bibitem{#2}}

\end{document}